\documentclass[12pt,a4paper]{amsart}
\usepackage{amsmath,amsfonts,amssymb,amscd,color,amsthm}
\usepackage[all]{xy}

\newtheorem{theorem}[equation]{Theorem}
\newtheorem{lemma}[equation]{Lemma}
\newtheorem{proposition}[equation]{Proposition}
\newtheorem{definition}[equation]{Definition}
\newtheorem{corollary}[equation]{Corollary}

\theoremstyle{remark}
\newtheorem{remark}[equation]{Remark}

\numberwithin{equation}{section}

\DeclareMathOperator{\Id}{Id}
\newcommand{\bbZ}{{\mathbb Z}}
\newcommand{\CS}{{\mathcal S}}
\newcommand{\uR}{{\underline R}}
\newcommand{\ur}{{\underline r}}
\newcommand{\ut}{{\underline t}}
\newcommand{\uw}{{\underline w}}
\newcommand{\uc}{{\underline c}}
\newcommand{\usigma}{{\underline \sigma}}
\newcommand{\bc}{{\bf c}}

\newcommand{\bs}{{\bf s}}
\newcommand{\bsigma}{{\boldsymbol\sigma}}
\newcommand{\inv}{^{-1}}

\newcommand{\eg}{{\it e.g.}}
\newcommand{\ie}{{\it i.e.}}
\newcommand{\lexp}[2]{\kern\scriptspace\vphantom{#2}^{#1}\kern-\scriptspace#2}
\newcommand{\Atilde}{G(\tilde A_{2n-1})}
\newcommand{\Mtilde}{M(\tilde A_{2n-1})}
\newcommand{\WA}{{W(\tilde A_{2n-1})}}
\newcommand{\WC}{{W(\tilde C_n)}}
\newcommand{\nnode}[1]{{\kern 5pt\hbox to
0pt{\hss{$\mathop\bigcirc\limits_{#1}$}\hss}\kern 5pt}}
\newcommand{\sbar}{{\vrule width10pt height3pt depth-2pt}}
\newcommand{\dbar}{{\rlap{\vrule width10pt height2pt depth-1pt} 
                 \vrule width10pt height4pt depth-3pt}}
\begin{document}
\author{F.~Digne}
\title[A Garside presentation for $\tilde C_n$]
{A Garside presentation for Artin-Tits groups of type $\tilde C_n$}
\address{LAMFA, CNRS et Universit\'e de Picardie-Jules Verne
33, Rue Saint-Leu 80039 Amiens Cedex France}
\email{digne@u-picardie.fr}
\urladdr{http://www.u-picardie.fr/digne}
\thanks{This work was partially supported by the ``Agence Nationale pour la Recherche''
project ``Th\'eories de Garside'' (number ANR-08-BLAN-0269-03)}

\begin{abstract}
We prove that an Artin-Tits group of type $\tilde C$ is the group of fractions of a Garside monoid, analogous
to the known dual monoids associated with Artin-Tits groups of spherical type and obtained by the ``generated group'' method.
This answers, in this particular case, a general question
on Artin-Tits groups, gives a new presentation of an Artin-Tits group of type $\tilde C$,
and has consequences for the word problem, the computation of some centralizers or the triviality of the center.
A key point of the proof is to show that this group is a group of fixed points in an Artin-Tits group of type
$\tilde A$ under an involution. Another important point is the study of the Hurwitz action of the usual
braid group on the decomposition of a Coxeter element into a product of reflections.
\end{abstract}
\maketitle
\section{Introduction}  The aim of this paper is to define a Garside structure on
the  Artin-Tits group  of type  $\tilde C_n$. A Garside  structure on a group
means  that this  group is  the group  of fractions  of a Garside monoid. A Garside monoid
is a monoid for which
the two posets given by right or
left divisibility have nice properties. 
In particular these two posets are lattices: there exist least common multiples and greatest common divisors,
and moreover these two lattices have a common sublattice
which generates the monoid and has a greatest element for both orders,
the Garside element (see Definition \ref{garside} below).
It is  known that all
Artin-Tits  groups of  spherical type  have two  nice Garside  structures
given respectively by the
classical monoid, obtained by generating the Artin-Tits group by lifts of the simple
reflections, and  the  dual  monoid  (see  \cite{bessis}), obtained by generating the
Artin-Tits group by elements lifting all reflections which divide (see below beginning
of Section \ref{tilde A}) a
given Coxeter element.  In  the  case  of
non-spherical  Artin-Tits groups  the classical  Artin-Tits monoid exists  but is only
locally  Garside (\ie, two elements have not always a common multiple,
in particular there is no Garside element).
An open  question in  general is  the existence  of a dual
Garside structure for general Artin-Tits groups. Such a structure is known for
type  $\tilde A$ and it has been  conjectured by  John Crisp
 and Jon McCammond  that no Artin-Tits group
of  affine type  other than  type $\tilde  A$ and  maybe $\tilde C$  can have such a
structure. A Garside structure provides normal forms for the elements
of the group and is a tool for solving the word problem. It also
allows to compute centralisers of periodic elements (roots of powers of the Garside element).

To get a dual Garside structure for an Artin-Tits group of type $\tilde C_n$,
we shall view this group  as the group
of  fixed points under  an involution in  an Artin-Tits group  of type $\tilde
A_{2n-1}$.  In \cite{Atilde} this  last  group has been shown to be
the group of fractions of several monoids 
only one of which, up to automorphism, is Garside, but unfortunately this one 
is not stable by the involution.
On the other hand only one of these non-Garside monoids, up to automorphism,
is stable by the involution. We show that
by taking fixed points in this last
monoid one gets a Garside structure for $\tilde C_n$.

The paper is organised as follows.
In Section \ref{germs} we introduce Garside monoids and give methods for getting a
Garside monoid from a partially defined product on a subset called a germ.
In  Section \ref{tilde A} we recall and improve  results from \cite{Atilde} on presentations of
Artin-Tits groups of type $\tilde A$. In  Section \ref{fixed points} we get a Garside structure
from  the fixed points of an involution in  a Garside group of type $\tilde A$.
In  Section \ref{tilde C_n} we show that an Artin-Tits group  of type $\tilde C$ can be seen as a
group of fixed points in an Artin-Tits group of type $\tilde A$.
In Section \ref{generated group} we show that the Garside structure we have got
in Section \ref{fixed points} can be obtained by the method of the ``generated group''
of \cite[0.4]{bessis}.
In Section \ref{hurwitz} we prove that the group of fractions of our Garside monoid
is the Artin-Tits group of type $\tilde C_n$ and we give a
a dual presentation of this group similar to what has been done in \cite{bessis},
\cite{freegroup} and \cite{Atilde} for
the other known dual monoids, where the generators are in one-to-one
correspondence with a set of reflections in the Coxeter group. One of the intermediate results
is that the
Hurwitz action is transitive on the set of shortest decompositions of a Coxeter element
of $\WC$ into a product of reflections. The analogous property is known
for all finite Coxeter groups (\cite[2.1.4]{bessis}), for all well-generated
complex reflexion groups (\cite[7.5]{Kpi1}), for Coxeter groups of type $\tilde A$
(\cite[3.4]{Atilde}) and is conjectured to be true for all Coxeter groups.
In section 8 we deduce from the Garside structure the centralizer of a power of
a lift of a Coxeter element in the Artin-Tits group.

I thank Eddy Godelle for having carefully read a preliminary version of this paper,
allowing me to fix an error in one of the proofs.

\section{Germs, Garside groups}\label{germs}
In this section we recall some definition and results on Garside monoids and groups.
\begin{definition}[\cite{catgar}]
\begin{enumerate}
\item
A germ of monoid is a set $P$ endowed with a partially defined product
$(x,y)\mapsto xy$, which
has a unit, \ie, an element 1 such that $1p$ and $p1$ are defined and equal
to $p$ for any $p\in P$.
\item A germ is associative if for any $a,b,c$ in $P$ such
that one of the products $a(bc)$ or $(ab)c$ is defined then the other one
is also defined and both products are equal.
\item A germ is left (resp.\ right) cancellative if $ab=ac$ (resp.\ $ba=ca$)
implies $b=c$. It is cancellative if it is cancellative on both sides.
\end{enumerate}
\end{definition}
We say that an element $a$ of a germ $P$ left divides an element $b\in P$
if there exists $c\in P$ such that $b=ac$.
Right divisibility is defined similarly.
In an associative germ right and left divisibility are preorder relations.
They are order relations if moreover the germ is cancellative and
there is no invertible element different from 1.
\begin{definition}
An associative germ is said to be left (resp.\ right) Noetherian
if there is no strictly decreasing infinite sequence
for left (resp.\ right) divisibility.
It is called Noetherian if it is both left and right
Noetherian.
\end{definition}
Note that
in a Noetherian germ there is no non-trivial invertible element.
Note also that to be left (resp.\ right) Noetherian is equivalent
to the fact that there is no strictly increasing bounded infinite sequence for
right (resp.\ left) divisibility.

A morphism from a germ to a monoid is a map which sends the (partial) product
on $P$ to the product in the monoid, and the unit of $P$ to the unit
of the monoid. The monoid $M(P)$ (resp.\ group $G(P)$) defined by a germ $P$ is
the monoid (resp.\ group) which has the universal property that it factorizes any
morphism from $P$ to a monoid (resp.\ group). In other words it is the monoid
(resp.\ group) generated by
$P$ with only relations the relations given by equalities of products in $P$.
It is known (\cite[3.5]{catgar})
that $P$ injects into $M(P)$ and 
is stable by left and right divisibility in $M(P)$.
The following definition is the definition of Garsideness that we will use in the
present paper (for small variations and generalizations of this definition see
\cite{dehornoy-paris}, \cite{dehornoy} and \cite{catgar}).

\begin{definition}\label{garside}
\begin{itemize}
\item
We say that a   monoid 
is Garside if it is cancellative, Noetherian, if it is a lattice
for both left and right divisibility
and if there exists an element $\Delta$ (called a Garside element)
whose sets of right and left divisors coincide and generate the monoid.
\item A group is  Garside  if it is generated by a submonoid which is a Garside monoid.
\end{itemize}
\end{definition}
Note that here we do not assume the set of divisors of $\Delta$ to be finite. When this set is infinite
what we call here a Garside monoid (resp.\ group) is what is usually called a quasi-Garside monoid
(resp.\  group).
A general reference for Garside
monoids can be \cite{catgar}.
The following result is a combination of \cite[3.31, 5.4 and 5.5]{catgar}.
\begin{proposition}\label{garside germ}
Let $P$ be an associative Noetherian germ satisfying the following properties:
\begin{enumerate}
\item two elements of $P$ have a least right
common multiple in $P$;
\item for all $m\in M(P)$, and $a$, $b$ in $P$, if $am=bm$ or
$ma=mb$, then $a=b$;
\item the elements of $P$ have a both left and right common multiple $\Delta\in P$.
\end{enumerate}
then $M(P)$ is a Garside monoid with Garside element $\Delta$.
\end{proposition}
Conversely, in a Garside monoid $M$ the divisors of the Garside element form a germ
$P$ satisfying the above properties and such that the canonical map $M(P)\to M$
is an isomorphism.
\begin{definition}
A germ satisfying the assumptions of \ref{garside germ} is called a Garside germ.
\end{definition}
Given a Garside germ $P$, elements in $M(P)$
have normal forms: any element can be written
uniquely $p_1p_2\ldots p_k$ with $p_i\in P$ such that $p_i$ is the greatest element for 
left divisibility dividing $p_ip_{i+1}\ldots p_k$. In the (Garside) group of a Garside
monoid $M$ with Garside element
$\Delta$ any element can be written as $\Delta^k x$ with $x\in M$ and $k\in \bbZ$.

We will use proposition \ref{garside germ} through its following corollary.
Before stating this corollary we need:
\begin{definition}\label{automorphisme de germe}
An automorphism of a germ $P$ is a bijection $f:P\to P$ mapping the unit to the
unit and such that $ab$ is
defined if and only if $f(a)f(b)$ is defined, in which case $f(ab)=f(a)f(b)$.
\end{definition}

\begin{corollary}\label{points fixes}
Let $P$ be an associative and Noetherian
germ having a (unique) both left and right common multiple $\Delta$;
assume that $M(P)$ is cancellative and that $P$ has
an automorphism $\sigma$ such that any two elements of $P$ have
a unique minimal $\sigma$-stable common right multiple; then
$P^\sigma$ is a Garside germ.
\end{corollary}
\begin{proof}
It is clear that $P^\sigma$ is an associative and Noetherian
germ. We get the result by 
proving that $P^\sigma$ satisfies the assumptions of Proposition
\ref{garside germ}.
Since $P$ is cancellative, if $x,y\in P^\sigma$ are such that $x$ divides $y$ in $P$ then
$x$ divides $y$ in $P^\sigma$. Hence the assumption of the corollary
implies that any two elements of $P^\sigma$
have an lcm in $P^\sigma$, whence (i).

The inclusion $P^\sigma\to M(P)$ extends
to a morphism $M(P^\sigma)\to M(P)$. If $am=bm$ or $ma=mb$ as in (ii), taking
the images in $M(P)$ we get $a=b$, since $M(P)$ is cancellative, whence (ii).

Unicity of $\Delta$ as the maximal element of $P$ implies $\Delta\in P^\sigma$, hence $\Delta$
is a common right and left multiple of $P^\sigma$,
whence (iii).
\end{proof}
\section{The monoids of type $\tilde A_{2n-1}$}\label{tilde A}
Before applying the previous
results to Artin-Tits groups of type
$\tilde C_n$, we need to recall
the presentation of an
Artin-Tits group of type $\tilde A$ given in \cite{Atilde}
as group of fractions of the
monoid $M(P)$ generated by an associative germ  $P$.

We recall first a general method for constructing a germ.
Given  a group $G$ generated as a monoid by a set $S$, \ie, any element
of $G$ is a product of elements of $S$, without inverses, we let
$l_S$ be the  length on $G$ with respect to $S$: the length of $g\in G$ is the
length of a shortest expression of $g$ as a product of elements of $S$. We
say  that $a\in G$ left divides $b\in G$, denoted by $a\preccurlyeq _Gb$
if $l_S(a)+l_S(a \inv b)=l_S(b)$, and
similarly for right divisibility. Starting with a balanced element $\delta$ (an
element  which has the same  set of right and  left divisors), we call $D$ the
set  of (left or right) divisors of $\delta$.  Then $D$ is a germ, the product
of  $a$  and  $b$  in  $D$  being  defined  and equal to $ab$ if $ab\in D$ and
$l_S(ab)=l_S(a)+l_S(b)$. Associated to this germ we have a monoid $M(D)$ and a
group  $G(D)$. If $D$ is a lattice then $M(D)$ is a Garside monoid
with Garside element $\delta$
(result due to J. Michel, see \cite[Theorem 0.5.2]{bessis}).
We call this construction the method of the ``generated group''.
This construction  starting with any  finite Coxeter group, its
set  of Coxeter generators and taking  the longest element for $\delta$, gives
the  associated Artin-Tits monoid  (or group). Starting  with a finite Coxeter
group  with set of generators all
reflections it gives the dual  monoid if  we take for  $\delta$ any Coxeter
element. Starting with a Coxeter element and all reflections in
a  Coxeter  group  of  type  $\tilde  A$ it gives the monoids $M(P)$ that we
describe in this section.
For  these  results see \cite{bessis},\cite{Atilde} and
\cite{catgar}.
  
We see a Coxeter group of type $\tilde A$
as a subgroup of the periodic permutations of $\bbZ$. We need some notation.
\begin{definition}
\begin{enumerate}\item
A permutation $w$ of $\bbZ$ is said to be $n$-periodic if
$w(i+n)=w(i)+n$ for any $i\in\bbZ$.
\item\label{cycle}
A cycle is an $n$-periodic
permutation which has precisely one orbit up to translation by
$n$. We say that a cycle is finite if its orbits are finite.
We call length of a cycle
the cardinality of one of its orbits.
We call support of
a permutation $w$ the union of its non-trivial orbits.
\end{enumerate}
\end{definition}

Any $n$-periodic permutation of $\bbZ$ can be written uniquely as a product of
disjoint cycles. There are two types of cycles: either all the orbits are
finite, or all the non-trivial orbits are infinite. In the
former case we will represent the cycle by one of its non-trivial orbits;
in the latter case we will represent a cycle as
$(a_1,a_2,\ldots,a_k)_{[h]}$, with all $a_i$ distinct modulo $n$,
meaning that the image of $a_i$ is $a_{i+1}$ for
$0\leq i<k$ and the image of $a_k$ is $a_1+nh$ (this cycle has $|h|$ non trivial orbits).
To each $n$-periodic permutation
$w$ of $\bbZ$ we can associate its total shift $\dfrac1n\sum_{x=1}^{x=n}(w(x)-x)$.
We recall the following facts:
\begin{proposition}\label{Coxeter A tilde}
The Coxeter group $W(\tilde A_{n-1})$ of type $\tilde A_{n-1}$
is isomorphic to 
the group of $n$-periodic permutations of $\bbZ$ with total shift equal to 0.
The reflections of $W(\tilde A_{n-1})$ correspond to
the permutations $(a,b)$ with $a$ and $b$ distinct modulo $n$.
The simple reflections are the permutations $s_i=(i,i+1)$ for $i=1,2,\ldots,n$. 
The reflections $s_i$ and $s_j$ commute unless $i-j=\pm1\pmod n$ in which case
their product has order 3.
\end{proposition}

In \cite{Atilde} a germ generating the Artin-Tits group of type $\tilde A_{n-1}$ is defined for each
partition of $\bbZ$ into two non-empty subsets $X$ and $\Xi$ stable by
translation by $n$; such a partition corresponds to the choice of a Coxeter element.
To recall this construction we need the following definitions.
A graphical representation of these definitions will be given after Definition \ref{pseudocycle}
We make the convention that Latin letters denote elements of
$X$ and Greek letters elements of $\Xi$.
\begin{definition}
We say that a cycle is positive and self non-crossing if it has one of
the following forms:
$(a_1,a_2,\ldots,a_k,\alpha_1,\alpha_2,\ldots,\alpha_l)$
or $(a_1,a_2,\ldots,a_k)_{[1]}$ or  $(\alpha_1,\alpha_2,\ldots,\alpha_l)_{[-1]}$
with $a_i\in X$, $\alpha_j\in \Xi$ satisfying the conditions
$a_1<a_2<\ldots<a_k<a_1+n$ and
$\alpha_1>\alpha_2>\ldots>\alpha_l>\alpha_1-n$.
\end{definition}
\begin{definition}
We say that two positive self non-crossing cycles are non-crossing if
they satisfy one of the following:
\begin{enumerate}
\item
One of them is of the form $(a_1,a_2,\ldots,a_k)$ or
$(a_1,a_2,\ldots,a_k)_{[1]}$ and the other one
$(\alpha_1,\alpha_2,\ldots,\alpha_l)$ or
$(\alpha_1,\alpha_2,\ldots,\alpha_l)_{[-1]}$. 
\item One of them is of the form
$(a_1,a_2,\ldots,a_k,\alpha_1,\alpha_2,\ldots,\alpha_l)$,
and the other $(b_1,\ldots,b_q,\beta_1,\ldots,\beta_r)$
with $k,l,q,r>0$,  $a_k<b_i<a_1+n$ for all
$i$ and  $\alpha_l<\beta_i<\alpha_1+n$ for all $i$.
\item One of them is of the form
$(a_1,a_2,\ldots,a_k,\alpha_1,\alpha_2,\ldots,\alpha_l)$
with $k,l\geq 0$,
and the other one  $(b_1,\ldots,b_m)$ (resp.\ $(\beta_1,\ldots,\beta_m)$)
with $a_j<b_i<a_{j+1}$ for some $j$ and all
$i$ (resp.\  $\alpha_{j+1}<\beta_i<\alpha_j$ for some $j$ and all $i$),
where in this condition we put $a_{k+1}=a_1+n$ and $\alpha_{l+1}=\alpha_1-n$,
(there is no condition on
$b_i$ if $k=0$ and no condition on $\beta_i$ if $l=0$).
\item One of them is of the form
$(a_1,a_2,\ldots,a_k)_{[1]}$ or $(\alpha_1,\alpha_2,\ldots,\alpha_k)_{[-1]}$
with $k>0$
and the other one  $(b_1,\ldots,b_m)$ (resp.\ $(\beta_1,\ldots,\beta_m)$)
with $a_j<b_i<a_{j+1}$ for some $j$ and all
$i$ (resp.\  $\alpha_{j+1}<\beta_i<\alpha_j$ for some $j$ and all $i$),
where in this condition we put $a_{k+1}=a_1+n$ and $\alpha_{k+1}=\alpha_1-n$.
\end{enumerate}
\end{definition}
\begin{definition}
We say that a periodic permutation of $\bbZ$ is positive and self non-crossing if it is the product
of disjoint positive self non-crossing and pairwise non-crossing cycles.
\end{definition}
Note that a finite cycle has total shift 0.
Note also that a periodic positive self non-crossing permutation of $\bbZ$
has at most 2 infinite cycles, one in $X$ and one in $\Xi$ and that, if it has total shift
0, it has either 0 or 2 infinite cycles.
\begin{definition}\label{pseudocycle}
We call pseudo-cycle a positive self non-crossing $n$-periodic
permutation of $\bbZ$ which
is the product of two infinite disjoint cycles.
\end{definition}
A pseudo-cycle is a permutation
$(a_1,\ldots,a_h)_{[1]}(\alpha_1,\alpha_2,\ldots,\alpha_l)_{[-1]}$
with the $a_i$ in $X$, the $\alpha_i$ in $\Xi$ and
$a_1<a_2<\ldots<a_h<a_1+n$,  $\alpha_1>\alpha_2>\ldots>\alpha_n>\alpha_1-n$.
A pseudo-cycle has total shift 0.

We  now give a graphical representation of the above definitions (see Figure 1
and  Figure 2; in these figures we have taken $n=9$, and modulo $n$ there are
5 elements in $X$ and $4$ in $\Xi$). 
Consider two  parallel oriented lines $D$  and $\Delta$ in the
plane,  with  same  orientation.  On  $D$  we  put a discrete set of points in
one-to-one  ordered correspondence with $X$ and  on $\Delta$ we put a discrete
set  of points in  one-to-one ordered correspondence  with $\Xi$. Let $\CS$ be
the  strip  delimited  in  the  plane  by  $D$ and $\Delta$. We associate to a
permutation  $w$ of $\bbZ$  a union of  oriented paths (one  for each orbit of
$w$)  contained in the  strip $\CS$ joining  $i$ to $w(i)$  for all $i$, up to
homotopy  in the  strip. A  cycle is  positive self  non-crossing if it can be
represented  by a  union of  disjoint oriented  paths such  that each  of them
intersects  $X$ (resp.\ $\Xi$) along an  increasing (resp.\ decreasing) or empty
subsequence.  Two positive self  non-crossing cycles are  non-crossing if they
can be represented by two unions of paths which do not cross each other.
Other figures (with $n=10$) can be found in Section \ref{generated group}.

\begin{figure}
\centerline{
\begin{xy}
0;<1cm,0cm>:
(1,0)="a3n" *={\bullet},*+!U{17},(2,0)="a2n" *={\bullet},*+!U{16},
(3,0)*={\bullet},*+!U{15},(4,0)="a1n" *={\bullet},*+!U{14},
(5,0)*={\bullet},*+!U{9},(6,0)="a3" *={\bullet},
*+!U{8},(7,0)="a2" *={\bullet},*+!U{7},
(8,0)*={\bullet},*+!U{6}, (9,0)="a1" *={\bullet}, *+!U{5},
(1,2)*={\bullet},*+!D{13},(2,2)="alpha1n" *={\bullet},*+!D{12},
(3,2)="alpha2n" *={\bullet},*+!D{11},
(4,2)*={\bullet},*+!D{10},
(6,2)*={\bullet},*+!D{4},
(7,2)="alpha1" *={\bullet},
*+!D{3},(8,2)="alpha2" *={\bullet},*+!D{2},
(9,2)*={\bullet},*+!D{1},
"alpha2";"a1"
**@{-};"a2"**\crv{(8,1)}?<>(.6)*\dir{>},"a2";"a3"**\crv{(6.5,1)};
"alpha1"**\crv{};"alpha2"**\crv{(7.5,1.7)}?<>(.6)*\dir{>},
"alpha2n";"a1n"
**@{-};"a2n"**\crv{(3,1)}?<>(.6)*\dir{>},
"a2n";"a3n"**\crv{(1.5,1)};"alpha1n"**\crv{};"alpha2n"
**\crv{(2.5,1.7)}?<>(.6)*\dir{>},
(10,2)*+!L{\Delta};(0,2) **\crv{}?>*\dir{>},
(10,0)*+!L{D};(0,0) **\crv{}?>*\dir{>}
\end{xy}
}
\caption{The $9$-periodic cycle $(5,7,8,3,2)$}
\end{figure}
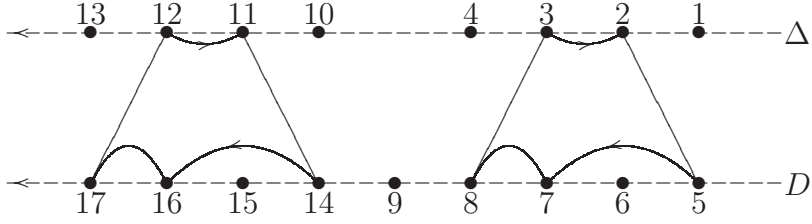

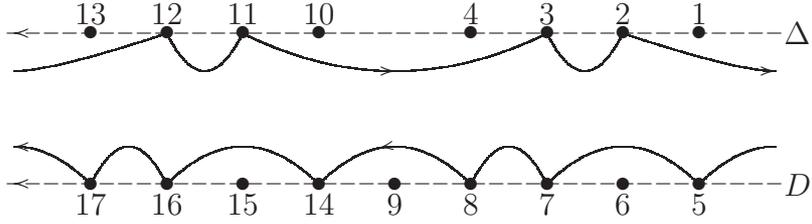
\begin{figure}
\centerline{
\begin{xy}
0;<1cm,0cm>:
(1,0)="a3n" *={\bullet},*+!U{17},(2,0)="a2n" *={\bullet},*+!U{16},
(3,0)*={\bullet},*+!U{15},(4,0)="a1n" *={\bullet},*+!U{14},
(5,0)*={\bullet},*+!U{9},(6,0)="a3" *={\bullet},
*+!U{8},(7,0)="a2" *={\bullet},*+!U{7},
(8,0)*={\bullet},*+!U{6}, (9,0)="a1" *={\bullet}, *+!U{5},
(1,2)*={\bullet},*+!D{13},(2,2)="alpha1n" *={\bullet},*+!D{12},
(3,2)="alpha2n" *={\bullet},*+!D{11},
(4,2)*={\bullet},*+!D{10},
(6,2)*={\bullet},*+!D{4},
(7,2)="alpha1" *={\bullet},
*+!D{3},(8,2)="alpha2" *={\bullet},*+!D{2},
(9,2)*={\bullet},*+!D{1},
(10,.5);"a1"**\crv{(9.5,.5)};
"a2"**\crv{(8,1)};"a3"**\crv{(6.5,1)};"a1n"**\crv{(5,1)}
?<>(.6)*\dir{>},"a1n";
"a2n"**\crv{(3,1)};"a3n"**\crv{(1.5,1)};(0,.5)**\crv{(.5,.5)}?>*\dir{>},
(10,1.5);"alpha2"**\crv{(9.5,1.5)}?<*\dir{<},"alpha2";"alpha1"**\crv{(7.5,1)};
"alpha2n"**\crv{(5,1)}?<>(.5)*\dir{<},
"alpha2n";"alpha1n"**\crv{(2.5,1)};(0,1.5)**\crv{(.5,1.5)},
(10,2)*+!L{\Delta};(0,2) **\crv{}?>*\dir{>},
(10,0)*+!L{D};(0,0) **\crv{}?>*\dir{>}
\end{xy}
}
\caption{The $9$-periodic pseudo-cycle
$(5,7,8)_{[1]}(3,2)_{[-1]}$}
\end{figure}

Following  \cite{Atilde} we now describe the  germ $P$ associated to the $X$ and $\Xi$.
The elements of $P$ are the $n$-periodic positive self non-crossing
permutations.
The permutation $c$ given by
$x_i\mapsto x_{i+1}$ and $\xi_i\mapsto \xi_{i-1}$ where $X=(x_i)_{i\in\bbZ}$
and $\Xi=(\xi_i)_{i\in\bbZ}$ with $x_i<x_{i+1}$ and $\xi_i<\xi_{i+1}$, is
a Coxeter element of $W(\tilde A_{n-1})$
and the elements  of $P$ are precisely the left (or right) divisors of $c$
(see \cite[Proposition 2.19]{Atilde}). The germ is obtained by the method of the generated group
described at the beginning of this section, hence the
product of two elements $w$ and  $w'$ of $P$ is defined in $P$ if
the  product $ww'$  of the  permutations is  in $P$ and if $l_{\tilde A_{n-1}}(ww')=l_{\tilde
A_{n-1}}(w)+l_{\tilde A_{n-1}}(w')$
where  $l_{\tilde A_{n-1}}$ is the length in $W(\tilde A_{n-1})$
with respect to the generating set consisting
of  all reflections. We will give below an equivalent and more tractable condition.

By general results $P$ is an associative, cancellative and
Noetherian germ, the monoid  $M(P)$ embeds in the group
$G(P)$ and any element $G(P)$ can be written $a\inv b$ with $a$ and $b$ in
$M(P)$.
Since the length $l_{\tilde A_{n-1}}$ is invariant by conjugation,
right and left divisibility in $P$ coincide.

The following is proved in \cite[Theorem 4.1]{Atilde}.
\begin{proposition}\label{G(P)=tilde A}
The group $G(P)$ is isomorphic to the Artin-Tits group of type $\tilde
A_{n-1}$.
\end{proposition}
Note that by \cite[Proposition 5.5]{Atilde} $P$ is a Garside germ if and only if either $X$ or
$\Xi$ contains exactly one element modulo $n$. 

We  now  give  a  more  tractable  definition  of divisibility in $P$. To each
element  $w\in P$  we associate  the partition  $p_w$ of $\bbZ$  whose parts are the
finite orbits of $w$ and the union of the two
infinite  orbits  of  $w$  if  they  exist.  Such  a partition is invariant by
translation  by $n$. We say that it is periodic. Moreover such a partition has
at most one infinite part and this infinite part must meet both $X$ and $\Xi$.
Also such a partition is non-crossing in the following sense:
\begin{definition}
\begin{itemize}
\item Two subsets $A$ and $B$ of  $X\cup\Xi$ are non-crossing  
if for any $a,a'\in A$ and any $b,b'\in B$ there exist in the strip 
$\CS$ a path $\gamma$ from $a$ to $a'$ and a path
$\delta$ from $b$ to $b'$ such that $\gamma$ and $\delta$ have an empty intersection.
\item We say that a partition is non-crossing if any two of its parts are non-crossing.
\end{itemize}
\end{definition}

By  \cite[Corollary 2.22]{Atilde}  $w\mapsto p_w$ is  a  bijection  from  $P$  onto  the  set  of
non-crossing  periodic partitions of $\bbZ$ such  that any infinite part meets
both  $X$ and $\Xi$. Such a partition can have at most one infinite part.
We order $P$ by
divisibility and the set of partitions by refinement. Note that the largest element of $P$
for the divisibility order is $c$ and that $p_c$ is the partition with only one part.

\begin{proposition}\label{p_w isomorphism}
The bijection $w\mapsto p_w$ is an isomorphism of ordered sets.
\end{proposition}
\begin{proof}  Let $v, w\in P$. We write  $v$ and $w$ as products $v=v_1\ldots
v_k$  and  $w=w_1\ldots  w_l$  where  each  $v_i$  and  each $w_i$ is either a
a positive self non-crossing finite cycle (see Definition \ref{cycle}) 
or is a pseudo-cycle (see Definition \ref{pseudocycle})
and the $v_i$ (resp.\ the $w_i$) are pairwise non-crossing.

By  \cite[Lemma 2.20]{Atilde} a reflection of $P$ divides $w$ in $W(\tilde A_n)$
if and only if its support is a subset of the
support  of  $w_i$  for  some  $i$.  Now  the following lemma
(see \cite[Corollary 2.9]{Atilde}) shows that
any  $v_i$  can be written as
a product $r_1r_2\ldots r_h$ in $P$ of
reflections whose union of supports is the support of $v_i$ and such
that the supports of $r_i$ and $r_{i+1}$ have a non-empty intersection.
\begin{lemma}\label{produit de reflexions}
The two following formulas give shortest decompositions of a finite cycle and
of a pseudo-cycle respectively into products of reflections of $\tilde A_{n-1}$.
\begin{multline*}
(a_1,a_2,\ldots,a_h)=(a_1,a_2)(a_2,a_3)\ldots (a_{h-1}a_h)\\
\shoveleft{(a_1,\ldots,a_h)_{[1]}(\alpha_1,\alpha_2,\ldots,\alpha_l)_{[-1]}=}\\
(a_1,a_2)(a_2,a_3)\ldots (a_{h-1},a_h)
(a_h,\alpha_1)(a_h,n+\alpha_1)
(\alpha_1,\alpha_2)(\alpha_2,\alpha_3) \ldots(\alpha_{l-1},\alpha_l)
\end{multline*}
\end{lemma}
\begin{proof}
By \cite[2.8]{Atilde} we know that
$l_{\tilde A_{n-1}}((a_1,\ldots,a_h)_{[1]}
(\alpha_1,\alpha_2,\ldots,\alpha_l)_{[-1]})=h+l$ and that
$l_{\tilde A_{n-1}}((a_1,a_2,\ldots,a_h))=h-1$.
\end{proof}
Assume  that $v=v_1\ldots v_k$ divides $w$ in $P$; fix $i$ and let us we write
$v_i=r_1\ldots r_h$ as above.
Then every $r_j$ divides $w$ so its support must
be  a subset of the support of some $w_s$. But since the supports of $r_j$ and
$r_{j+1}$  have  a  non-empty  intersection,  they  have to be included in the
support  of the same $w_s$, so that the  whole support of $v_i$ is a subset of
the  support of some $w_s$. This, being true  for all $i$, means that $p_v$ is
finer than $p_w$.

Conversely,  assume that $p_v$ is  finer than $p_w$. We  prove by induction on
$l_{\tilde A_{n-1}}(v)$ that $v$  divides $w$ in $P$.  If $v$ is trivial the result is true.
If $v$ is not trivial, there
exists  a reflection $s$ dividing $v$, so  that the support of $s$ is included
in the support of some $v_i$. Hence the support of $s$ is also included in the
support  of some $w_j$, so that $s$ divides $w$. Put $v=sv'$ and $w=sw'$; then
$v$ divides $w$ if and only if $v'$ divides $w'$. We will be done by induction
if  we  prove  that  $p_{v'}$  is  finer  than $p_{w'}$. We have seen that $s$
divides  a cycle or  a pseudo-cycle $v_i$  of $v$ and  a cycle or pseudo-cycle
$w_j$  of $w$. Since disjoint cycles and pseudo-cycles commute we may assume that
$s$ divides $v_1$ and $w_1$. We have only to
show   that  $p_{v'_1}$  is  finer   than  $p_{w'_1}$  where  $v_1=sv'_1$  and
$w_1=sw'_1$. This will be a consequence of the following lemma:
\begin{lemma} Let $u\in P$ be a finite cycle or a pseudo-cycle and let $s$ be
a reflection dividing $u$ in $P$. We put $u=su'$.
\begin{itemize}
\item If $u=(a_1,a_2,\ldots,a_h)\in P$ is a finite cycle (with $a_i\in X\cup \Xi$)
and $s=(a_1,a_j)$, then $u'=(a_1,a_2,\ldots,a_{j-1})(a_j,a_{j+1},\ldots,a_h)$.
\item  If $u=(a_1,\ldots,a_h)_{[1]}(\alpha_1,\alpha_2,\ldots,\alpha_l)_{[-1]}$ with $a_i\in X$ 
and $\alpha_i\in\Xi$ is a pseudo-cycle and $s=(a_1,a_j)$ then
$u'=(a_1,a_2,\ldots,a_{j-1})
(a_j,a_{j+1},\ldots,a_h)_{[1]}(\alpha_1,\alpha_2,\ldots,\alpha_l)_{[-1]}$.
\item  If $u=(a_1,\ldots,a_h)_{[1]}(\alpha_1,\alpha_2,\ldots,\alpha_l)_{[-1]}$ with $a_i\in X$ 
and $\alpha_i\in\Xi$ is a pseudo-cycle and $s=(\alpha_1,\alpha_j)$ then
$u'=(a_1,a_2,\ldots,\ldots,a_h)_{[1]}(\alpha_1,\alpha_2,\ldots,\alpha_{j-1})
(\alpha_j,\ldots,\alpha_l)_{[-1]}$.
\item  If $u=(a_1,\ldots,a_h)_{[1]}(\alpha_1,\alpha_2,\ldots,\alpha_l)_{[-1]}$ with $a_i\in X$ 
and $\alpha_i\in\Xi$ is a pseudo-cycle and $s=(a_1,\alpha_1)$ then
$u'=(a_1,a_2,\ldots,\ldots,a_h,\alpha_1+n,\alpha_2+n,\ldots,\alpha_l+n)$.
\end{itemize}
\end{lemma}
This lemma is an easy computation (see also \cite[Lemma 2.5]{Atilde}).

Applying the lemma with $u=v_1$ and with $u=w_1$ shows that $p_{v'_1}$ is finer than $p_{w'_1}$.
\end{proof}

\section{Fixed points in $\tilde A_{2n-1}$}\label{fixed points}

It  is well known that the Coxeter group  of type $\tilde C_n$ is the group of
fixed  points  under  the  involution  of  $\WA$ which maps $s_i$ to $s_{2n-i}$,
with the notation of Proposition \ref{Coxeter A tilde},
the subscript $i$ being taken modulo $2n$ (see \eg, \cite{hee}).
This involution can be lifted to the Artin-Tits group using the same formula.
We shall see in Section \ref{tilde C_n} that
similarly the Artin-Tits group of type $\tilde C_n$ is the group of fixed
points under this lifted involution in the Artin-Tits group of type $\tilde A_{2n-1}$.

To get a Garside germ for an Artin-Tits group of type $\tilde C_n$ we shall
start with a particular choice of $X$ and $\Xi$ in the construction of the
previous section, compatible with this involution. We take $X$ to be the the set of
odd integers and $\Xi$ to be the set of
even integers. This corresponds to choosing the Coxeter element $c=
(2,3)(4,5)\ldots (2n,2n+1)(1,2)(3,4)\ldots(2n-1,2n)=
s_2s_4\ldots s_{2n}s_1s_3\ldots s_{2n-1}$ in the previous section.
Then the germ $P$ has an automorphism $\sigma$ coming from the
map $i\mapsto 1-i$ which interchanges $X$ and $\Xi$. This involution lifts to the Artin-Tits group the
involution $s_i\mapsto s_{2n-i}$ of the Coxeter group. We still denote by $\sigma$
this involution of $X\cup\Xi$.
\begin{theorem}
The germ $P^\sigma$ is Garside
\end{theorem}
\begin{proof}
We show that $P$ and $\sigma$ satisfy the assumptions of Corollary
\ref{points fixes}.

First $P$ has a unique both right and left multiple of all
its elements, which is $c$. Note that
with our choice of $X$ and $\Xi$, the element $c$ seen as a $2n$-periodic
permutation of $\bbZ$ is
$$i\mapsto \begin{cases}i+2&\text{ for odd } i,\\
i-2&\text{ for even } i.\end{cases}$$
We have to show that any two elements $v$ and $w$ in $P$ have a unique minimal
$\sigma$-stable common multiple $z$. By Proposition \ref{p_w isomorphism}, this
amounts to show the existence of a unique minimal  $\sigma$-stable non-crossing
periodic partition of the form $p_z$ coarser than $p_v$ and $p_w$. Note that a 
$\sigma$-stable non-crossing partition has precisely 0, one or two infinite parts.
Hence the condition for such a partition to be of the form $p_z$, \ie,
that any infinite part has a non empty
intersection with both $X$ and $\Xi$ is equivalent to the condition for the partition
to have at most one infinite part.

A  $\sigma$-stable partition coarser than $p_v$ and $p_w$ is also coarser than
the  partitions $\lexp\sigma p_v$  and $\lexp\sigma p_w$.  We claim that it is
sufficient  to show the  existence of a  unique minimal non-crossing partition
coarser  than  these  4  partitions:  such  a  partition  will be periodic and
$\sigma$-stable by unicity. Moreover since it is non-crossing and symmetric it
will  have precisely  0, one  or two  infinite parts.  If it  has zero  or one
infinite  part it is of the form $p_z$ and we are done. If it has two infinite
parts these parts are interchanged by $\sigma$, in which case there is exactly
one  minimal coarser  partition with  one infinite  part meeting  both $X$ and
$\Xi$,  obtained by putting together the two infinite parts. This partition is
$\sigma$-stable   and  periodic.  By   Proposition  \ref{p_w  isomorphism}  it
corresponds to a minimal $\sigma$-stable common multiple of $v$ and $w$ and we
are done also in this case.

The theorem is thus a consequence of the following lemma which ensures by induction
the existence 
of a unique non-crossing partition coarser than any (finite) given number of
partitions of $X\cup\Xi$.
\begin{lemma}
There exists a unique minimal non-crossing partition coarser than any two given partitions
of $X\cup \Xi$.
\end{lemma}
\begin{proof}[Proof of the lemma]
Consider a graph whose vertices are the parts of the given partitions and such that
there is an edge between two vertices if the corresponding parts cross each other
(this includes the case where two parts have an intersection). Then the desired partition
is the partition whose parts are the union of all the parts lying in the same connected component of our graph.
\end{proof}
\end{proof}
\section{The Artin-Tits group of type $\tilde C_n$}\label{tilde C_n}\label{Artin-Tits de type tilde C}
We denote by $\Atilde$ the Artin-Tits group of type $\tilde A_{2n-1}$ and
by $G(\tilde C_n)$ the Artin-Tits group of type $\tilde C_n$. They are the
groups of fractions of the corresponding classical Artin-Tits monoids, respectively
$\Mtilde$ and $M(\tilde C_n)$. The monoid
$\Mtilde$ has a presentation with generators
$\bs_1,\ldots,\bs_{2n}$ and relations $\bs_is_j=\bs_j\bs_i$ if $|i- j|\neq 1 \text{ mod } 2n$
and $\bs_i\bs_{i+1}\bs_i=\bs_{i+1}\bs_i\bs_{i+1}$ for all $i$,
where the indices are taken modulo $2n$; the monoid 
$M(\tilde C_n)$  of type $\tilde C_n$ has a presentation with generators
$\bsigma_0,\ldots,\bsigma_n$
and relations the braid relations given by the Coxeter diagram
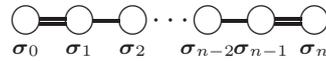
\begin{equation*}
\nnode
{\bsigma_0}\dbar\nnode{\bsigma_1}\sbar\nnode{\bsigma_2}
\cdots\nnode{\bsigma_{n-2}}\sbar\nnode{\bsigma_{n-1}}\dbar\nnode{\bsigma_n}\kern 4pt
\end{equation*}
Like any Artin-Tits monoid, by \cite{paris}, the monoids
$\Mtilde$ and $M(\tilde C_n)$ embed in their respective
groups $\Atilde$ and $G(\tilde C_n)$.

By Proposition \ref{G(P)=tilde A} the group $\Atilde$ is
isomorphic to the group of fractions of the monoid $M(P)$ 
generated by the germ $P$ defined in Section \ref{fixed points}. 
This isomorphism maps the generator $\bs_i$ of $\Mtilde$ to the element $(i,i+1)\in P$,
so that $\Mtilde$ is a submonoid of $M(P)$. The involution $\sigma$ considered in Section
\ref{fixed points} restricts to $\Mtilde$ in the diagram automorphism of $\tilde A_{2n-1}$
which maps $\bs_i$ to $\bs_{2n-i}$ where the indices are taken modulo $2n$.

By \cite[4.4]{michel} the monoid of fixed points $\Mtilde^\sigma$ 
is isomorphic to $M(\tilde C_n)$. We will identify these two monoids. Under this
identification we have
$\bsigma_0=\bs_{2n}$, $\bsigma_1=\bs_1\bs_{2n-1}$,
\dots, $\bsigma_{n-1}=\bs_{n-1}\bs_{n+1}$, $\bsigma_n=\bs_n$.

We summarise all these facts in the commutative diagram:
\begin{equation}\label{diagramme}
\xymatrix{
M(\tilde
C_n)\ar[r]^\sim\ar@{^{(}->}[dd]&\Mtilde^\sigma\ar@{^{(}->}[r]\ar@{^{(}->}[d]&\Mtilde\ar@{^{(}->}[d]\\
&M(P)^\sigma\ar@{^{(}->}[r]\ar@{^{(}->}[d]&M(P)\ar@{^{(}->}[d]\\
G(\tilde C_n)\ar[r]&\Atilde^\sigma\ar@{^{(}->}[r]&\Atilde=G(P)
}
\end{equation}
The aim of this section is to prove the following theorem.
I thank John Crisp and Dave Margalit for having pointed out that this theorem can be easily
deduced from the results of Birman-Hilden and
Maclachlan-Harvey.
\begin{theorem}\label{G(tilde C)=C_Atilde(sigma)}
The morphism $G(\tilde C_n)\rightarrow\Atilde^\sigma$ of Diagram \ref{diagramme}, which maps 
$\bsigma_0$ to $\bs_{2n}$, $\bsigma_i$ to $\bs_i\bs_{2n-i}$ for $i=1,\ldots,n-1$
and $\bsigma_n$ to $\bs_n$ is an isomorphism.
\end{theorem}
We need first to recall how the Artin-Tits groups of type $\tilde A$ and $\tilde C$ can be embedded
into mapping class
groups. Let $E$ and $F$ be two finite subsets of a 2-sphere $S^2$.
We denote by $M(S^2,E,F)$ the subgroup of the mapping class group of the 2-sphere
with punctures the points of $E\cup F$, represented
by the diffeomorphisms which fix each point of $F$ and stabilize $E$.
The following result appears in \cite[Section 2]{charney-crisp}:
\begin{proposition}[Charney-Crisp]
\begin{enumerate}
\item\label{charney-crisp1} Let $E=\{P_1,\ldots, P_n\}$ and
$F=\{P',P''\})$ be two sets of points of $S^2$; the group $G(\tilde A_{n-1})$
embeds into the group $M(S^2,E,F)$. This embedding maps the standard
generator $\bs_i$ of $G(\tilde A_{n-1})$ to the positive braid twist exchanging $P_i$ and $P_{i+1}$,
for $1\leq i\leq n$, where the indices are taken modulo $n$.
\item\label{charney-crisp2} Let $E=\{P_1,\ldots,P_n\}$ and $F=\{P_0,Q',Q''\}$ be subsets of $S^2$
then the group $G(\tilde C_n)$ is isomorphic to $M(S^2,E,F)$.
This isomorphism maps the standard
generator (numbered as in the beginning of Section \ref{tilde C_n})
$\bsigma_i$ of $G(\tilde C_n)$ to the positive braid twist exchanging $P_i$ and $P_{i+1}$ for
$1\leq i\leq n-1$ and maps $\bsigma_0$ (resp.\ $\bsigma_n$) to the square
of the positive braid twist exchanging
$P_0$ and $P_1$ (resp.\ $P_n$ and $Q''$).
\end{enumerate}
\end{proposition}
We need also the following result of Birman and Hilden (see \cite{birman-hilden} and \cite{harvey-mac}).
\begin{theorem}[Birman-Hilden]\label{birman-hilden}
Let $E_1$ and $E_2$ be two finite sets of points of a 2-sphere $S^2$, let
$\pi: S^2\to S^2/G$ be a ramified covering realizing the quotient of $S^2$ by a finite group $G$
of diffeomorphisms stabilizing $E_1$ and $E_2$
and let $F\subset S^2/G$ be the set of ramification points. Then the projection induces
an isomorphism from the normalizer of $G$ in the mapping class
group $M(S^2,E_1,E_2)$ to the mapping class group $M(S^2/G,\pi(E_1),\pi(E_2)\cup F)$.
\end{theorem}
\begin{proof}[Proof of Theorem \ref{G(tilde C)=C_Atilde(sigma)}]
Consider a set $E$ consisting of
$2n$ points $P_1,\ldots,P_{2n}$ regularly placed on the
equator of a sphere $S^2$ and let $P'$ and $P''$ be the north and south poles.
We apply Proposition \ref{charney-crisp1}, for embedding
the group $G(\tilde A_{2n-1})$ into the mapping class group $M(S^2,E,F)$ where
$F=\{P',P''\}$.
We then apply Theorem \ref{birman-hilden} with $G$ the group of order 2 generated by the symmetry $\sigma$ which exchanges $P'$
and $P''$ and exchanges $P_i$ and $P_{2n-i}$ for all $i$. We can view $S^2/G$ as a sphere so that
$\pi:S^2\to S^2/G$ is a ramified covering of a 2-sphere
with 2 ramification points  $A',A''\in S^2/G$ which are antipodal on the equator. The points
$\pi(P_i)=\pi(P_{2n-i})$
are regularly placed on one half of the equator and the point $P=\pi(P')=\pi(P'')$ is
on the other half of this equator. We get an isomorphism 
$M(S^2,E,F)^\sigma\simeq M(S^2/G,\pi(E),\{P,A',A''\})$. Since $\pi(E)$ has cardinality
$n$, we get by Proposition \ref{charney-crisp2} that  $M(S^2/G,\pi(E),\{P,A',A''\})$
is isomorphic to $G(\tilde C_n)$. Since $\sigma$ maps the product of braid twists
$\bs_i\bs_{n-i}$ to the braid twist
$\bsigma_i$ for $i=1,\ldots,n-1$ and maps the braid twist $\bs_n$ (resp.\ $\bs_{2n}$) to the
square $\bsigma_0$ (resp.\ $\bsigma_n$) of the positive braid twist exchanging $A'$ and
$\pi(P_1)$ (resp.\ $\pi(P_n)$ and $A''$), we get that the isomorphism restricted to $\Atilde^\sigma$
is onto which means that $\Atilde^\sigma$ is in fact equal to $M(S^2,E,F)^\sigma$.
Moreover we also see that the above isomorphism coincides with the
morphism of Theorem \ref{G(tilde C)=C_Atilde(sigma)}, whence the result.
\end{proof}

\section{A generated group}\label{generated group}
Our aim in  this section is  to show that the germ $P^\sigma$ is obtained from
the Coxeter group of type $\tilde C_n$ by the method of the generated group
described at the beginning of Section \ref{tilde A}.
We  denote the  Coxeter group  of type  $\tilde C_n$ by
$\WC$  and see it as the  group of  fixed points  under the involution
$\sigma:s_i\mapsto
s_{2n-i}$  in  the  Coxeter  group  $\WA$  of  type $\tilde A_{2n-1}$ (see the
beginning of Section \ref{fixed points}).
\begin{lemma}\label{reflexions de tilde C}
The reflections of $\WC$ are the $\sigma$-stable reflections of $\WA$ and 
the products $r.\lexp\sigma r$ where $r$ is a reflection of $\WA$ which is
not $\sigma$-stable.
\end{lemma}
\begin{proof}
This  is well known. We recall a  proof. By \cite{hee},
if $\sigma$ is an automorphism  of a Coxeter group 
permuting the simple reflections, the group
of  $\sigma$-fixed elements  is a  Coxeter group  with simple  reflections the
longest  elements of the parabolic  subgroups generated by the $\sigma$-orbits
of spherical type of simple reflections. The simple reflections of $\WC$ are thus $s_n$, $s_{2n}$ and
the  products  $s_i.\lexp\sigma  s_i=s_is_{2n-i}$  for  $1\leq i\leq n-1$.
An arbitrary reflection of $\WC$ is a conjugate under $\WC$ to a
simple  reflection,  whence  the  result  since this conjugation commutes with
$\sigma$.
\end{proof}
We  put $c=s_2s_4\ldots s_{2n}s_1s_3\ldots s_{2n-1}$  as in Section \ref{fixed 
points}. We denote the simple reflections of $\WC$ by $\sigma_0=s_{2n}$,
$\sigma_1=s_1s_{2n-1}$, $\sigma_2=s_2s_{2n-2}$,\dots,
$\sigma_{n-1}=s_{n-1}s_{n+1}$ and $\sigma_n=s_n$. We have
$c=\sigma_0\sigma_2\sigma_4\ldots\sigma_1\sigma_3\sigma_5\ldots$, hence $c$ is
both a $\sigma$-fixed Coxeter element of $\WA$ and a Coxeter element of $\WC$.
We  denote by $\preccurlyeq_{\WA}$ (resp.\ $\preccurlyeq_{\WC}$) the
divisibility  in $\WA$ (resp.\  $\WC$) with respect  to the reflection length.
Beware that the reflection length in $\WC$ is not the restriction of the
reflection length in $\WA$, in particular the reflection length of $c$ is
$n+1$ in $\WC$ and $2n$ in $\WA$, since a Coxeter element in a Coxeter group
of rank $h$ has reflection length equal to $h$ (see \eg, \cite{dyer}).
Recall that divisibility  for the reflection length in a Coxeter group
can be  defined by saying that $w$ divides
$w'$ if for one, or for any, shortest decomposition of $w$ into a product of
reflections, there exists a shortest decomposition of $w'$ as a product of
reflections which begins with that decomposition
of $w$. We now show  that divisibility
in $\WA$ restricts to divisibility in $\WC$ for the divisors of $c$.
Recall that by the results of \cite{Atilde} $P$ is the set of divisors of $c$ in
$\WA$, so that $P^\sigma$ is the set of elements of $\WC$ which divide $c$ in
$\WA$.

\begin{proposition}\label{divisibility restricts}
Let $w$ and $w'$ be two elements of $\WC$; then 
$w\preccurlyeq_{\WC}w'\preccurlyeq_{\WC}c$ if and only if
$w\preccurlyeq_{\WA} w'\preccurlyeq_{\WA} c$.
\end{proposition}
\begin{proof}
We first prove that divisibility in $\WC$ implies divisibility in $\WA$ for divisors of $c$
in $\WC$.

We claim that in any decomposition of $c$ into a product of reflections of $\WC$
there is at least one reflection conjugate to $\sigma_0$ and one reflection conjugate to
$\sigma_n$. Indeed the quotient of $\WC$ by its commutator subgroup is the direct product
of 3 groups of order 2, generated respectively by the image of $\sigma_0$, the image of
$\sigma_n$ and the image of any one of
$\sigma_1$, \dots, $\sigma_{n-1}$ which have all same image since they are conjugate.
Hence the image of $c$ in this quotient must involve the
images of $\sigma_0$ and $\sigma_n$. Whence the result since two conjugate reflections have
same image in this quotient.

Let $c=r_0\ldots r_n$ be a decomposition of $c$ into a product of $n+1$
reflections in $\WC$. Let $k$ be the
number of reflections in that decomposition which are conjugate to either
$\sigma_0$ or $\sigma_n$.
These $k$ reflections are also reflections in $\WA$ and the other $n+1-k$ reflections are products
of 2 reflections in $\WA$. So this product is the product of $k+2(n+1-k)=2n+2-k$ reflections in $\WA$. But
$c$ cannot be written with less than $2n$ reflections in $\WA$, whence $k\leq 2$. Since we know that
$k\geq 2$, we have $k=2$. We have proved that, replacing a reflection of $\WC$ by its image in
the embedding $\WC\hookrightarrow\WA$,
any shortest decomposition of $c$ as a product of reflection of $\WC$
becomes a shortest decomposition of $c$ as a product of reflections of $\WA$.
This gives that if
$w\preccurlyeq_{\WC}w'\preccurlyeq_{\WC}c$ then
$w\preccurlyeq_{\WA} w'\preccurlyeq_{\WA} c$.

To prove the converse we first prove
a formula for the reflection length $l_\WC$ in $\WC$ of the elements of $P^\sigma$.
An element of $P^\sigma$ can be written uniquely as a commuting
product of elements of three types:
either $\sigma$-stable non-crossing finite cycles, or products of a  finite cycle
and its image by $\sigma$, or $\sigma$-stable pseudo-cycles.
There can be at most
one pseudo-cycle and any two of these factors are non-crossing.
In Figures 3, 4 and 5 we have given examples
in $\tilde C_5$, \ie, the partitions are $10$-periodic
and $\sigma$ comes from $i\mapsto 1-i \text{ modulo }10$.
\begin{figure}
\centerline{
\begin{xy}
0;<1cm,0cm>:(1,0)*={\bullet},*+!U{19},
(2,0) *={\bullet},*+!U{17},(3,0)="a3n" *={\bullet},*+!U{15},
(4,0)="a2n" *={\bullet},*+!U{13},(5,0)="a1n" *={\bullet},*+!U{11},
(6,0)*={\bullet},*+!U{9},(7,0) *={\bullet},
*+!U{7},(8,0)="a3" *={\bullet},*+!U{5},
(9,0)="a2" *={\bullet},*+!U{3}, (10,0)="a1" *={\bullet}, *+!U{1},
(1,2)="alpha1n" *={\bullet},*+!D{20},
(2,2)="alpha2n" *={\bullet},*+!D{18},(3,2)="alpha3n" *={\bullet},*+!D{16},
(4,2) *={\bullet},*+!D{14},
(5,2)*={\bullet},*+!D{12},
(6,2)="alpha1" *={\bullet},*+!D{10},
(7,2)="alpha2" *={\bullet},*+!D{8},
(8,2)="alpha3" *={\bullet},
*+!D{6},(9,2) *={\bullet},*+!D{4},
(10,2)*={\bullet},*+!D{2},
"alpha3";"a1"
**@{-};"a2"**\crv{(9.5,.5)}?<>(.6)*\dir{>},"a2";"a3"**\crv{(8.5,.5)};
"alpha1"**\crv{};"alpha2"**\crv{(6.5,1.7)}?<>(.6)*\dir{>};"alpha2";
"alpha3"**\crv{(7.5,1.7)}?<>(.6)*\dir{>},
"alpha3n";"a1n"
**@{-};"a2n"**\crv{(4.5,.5)}?<>(.6)*\dir{>},"a2n";"a3n"**\crv{(3.5,.5)};
"alpha1n"**\crv{};"alpha2n"**\crv{(1.5,1.7)}?<>(.6)*\dir{>};"alpha2n";
"alpha3n"**\crv{(2.5,1.7)}?<>(.6)*\dir{>},
(11,2)*+!L{\Delta};(0,2) **\crv{}?>*\dir{>},
(11,0)*+!L{D};(0,0) **\crv{}?>*\dir{>}
\end{xy}
}
\caption{A $\sigma$-stable cycle in $\tilde A_9$}
\end{figure}
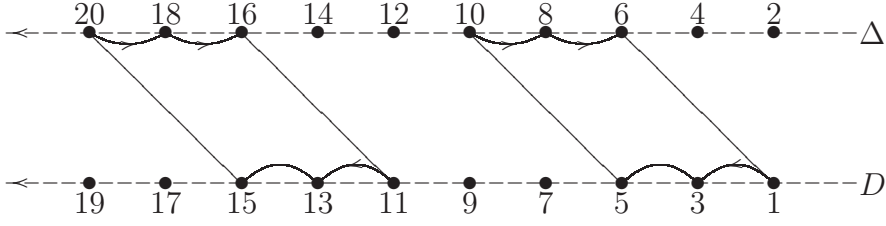

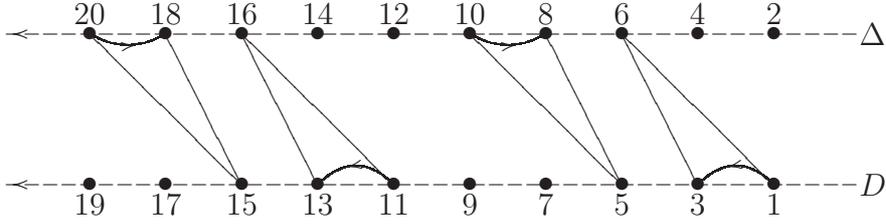
\begin{figure}
\centerline{
\begin{xy}
0;<1cm,0cm>:(1,0)*={\bullet},*+!U{19},
(2,0) *={\bullet},*+!U{17},(3,0)="a3n" *={\bullet},*+!U{15},
(4,0)="a2n" *={\bullet},*+!U{13},(5,0)="a1n" *={\bullet},*+!U{11},
(6,0)*={\bullet},*+!U{9},(7,0) *={\bullet},
*+!U{7},(8,0)="a3" *={\bullet},*+!U{5},
(9,0)="a2" *={\bullet},*+!U{3}, (10,0)="a1" *={\bullet}, *+!U{1},
(1,2)="alpha1n" *={\bullet},*+!D{20},
(2,2)="alpha2n" *={\bullet},*+!D{18},(3,2)="alpha3n" *={\bullet},*+!D{16},
(4,2) *={\bullet},*+!D{14},
(5,2)*={\bullet},*+!D{12},
(6,2)="alpha1" *={\bullet},*+!D{10},
(7,2)="alpha2" *={\bullet},*+!D{8},
(8,2)="alpha3" *={\bullet},
*+!D{6},(9,2) *={\bullet},*+!D{4},
(10,2)*={\bullet},*+!D{2},
"alpha3";"a1"
**@{-};"a2"**\crv{(9.5,.5)}?<>(.6)*\dir{>},"a2";"alpha3"**@{-},
"alpha3n";"a1n"
**@{-};"a2n"**\crv{(4.5,.5)}?<>(.6)*\dir{>},"a2n";"alpha3n"**@{-},
"a3";"alpha1"
**@{-};"alpha2"**\crv{(6.5,1.7)}?<>(.6)*\dir{>},"alpha2";"a3"**@{-},
"a3n";"alpha1n"
**@{-};"alpha2n"**\crv{(1.5,1.7)}?<>(.6)*\dir{>},"alpha2n";"a3n"**@{-},
(11,2)*+!L{\Delta};(0,2) **\crv{}?>*\dir{>},
(11,0)*+!L{D};(0,0) **\crv{}?>*\dir{>}
\end{xy}
}
\caption{A non-crossing product of a cycle and its image by $\sigma$ in $\tilde A_9$}
\end{figure}

\begin{figure}
\centerline{
\begin{xy}
0;<1cm,0cm>:(1,0)*={\bullet},*+!U{19},
(2,0)="a3n" *={\bullet},*+!U{17},(3,0) *={\bullet},*+!U{15},
(4,0)="a2n" *={\bullet},*+!U{13},(5,0)="a1n" *={\bullet},*+!U{11},
(6,0)*={\bullet},*+!U{9},(7,0)="a3" *={\bullet},
*+!U{7},(8,0) *={\bullet},*+!U{5},
(9,0)="a2" *={\bullet},*+!U{3}, (10,0)="a1" *={\bullet}, *+!U{1},
(1,2)="alpha1n" *={\bullet},*+!D{20},
(2,2)="alpha2n" *={\bullet},*+!D{18},(3,2) *={\bullet},*+!D{16},
(4,2)="alpha3n" *={\bullet},*+!D{14},
(5,2)*={\bullet},*+!D{12},
(6,2)="alpha1" *={\bullet},*+!D{10},
(7,2)="alpha2" *={\bullet},*+!D{8},
(8,2) *={\bullet},
*+!D{6},(9,2)="alpha3" *={\bullet},*+!D{4},
(10,2) *={\bullet},*+!D{2},
"a1";"a2"**\crv{(9.5,.5)}?<>(.6)*\dir{>},"a2";"a3"**\crv{(8,.5)};"a3";"a1n"**\crv{(6,.5)},
"alpha1";"alpha2"**\crv{(6.5,1.7)}?<>(.6)*\dir{>};"alpha2";
"alpha3"**\crv{(8,1.7)}?<>(.6)*\dir{>},
"a1n";"a2n"**\crv{(4.5,.5)}?<>(.6)*\dir{>},"a2n";"a3n"**\crv{(3,.5)};
"alpha1n";"alpha2n"**\crv{(1.5,1.7)}?<>(.6)*\dir{>};"alpha2n";
"alpha3n"**\crv{(3,1.7)}?<>(.6)*\dir{>},"alpha3n";"alpha1"**\crv{(5,1.7)},
"alpha1n";(0,1.7)**\crv{(.5,1.7)},
"alpha3";(11,1.7)**\crv{(10.5,1.7)},
"a3n";(0,.5)**\crv{(.5,.5)},
"a1";(11,.5)**\crv{(10.5,.5)},
(11,2)*+!L{\Delta};(0,2) **\crv{}?>*\dir{>},
(11,0)*+!L{D};(0,0) **\crv{}?>*\dir{>}
\end{xy}
}
\caption{A $\sigma$-stable pseudo-cycle in $\tilde A_9$}
\end{figure}
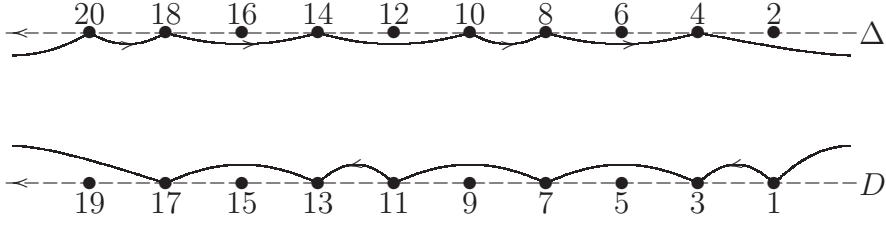

\begin{lemma}\label{formule pour l dans Ctilde}
Let $w=w_1\ldots w_r$ be a decomposition as above of $w\in P^\sigma$; then $l_\WC(w)=\sum_{i=1}^{i=r}l_\WC(w_i)$
and $$l_\WC(w_i)=
\begin{cases}
k&\text{ if $w_i$ is a $\sigma$-stable cycle of length $2k$}\\
k-1&\text{ if $w_i$ is the product of a cycle of length $k$ and its image by $\sigma$}\\
k+1&\text{ if $w_i$ is a pseudo-cycle whose support modulo $2n$ has cardinality $2k$}
\end{cases}
$$
\end{lemma}
\begin{proof}
We denote by $f(w)$ the formula of the lemma. We have to prove that $l_\WC(w)=f(w)$.
First, it is clear that $f(w)=0$ if and only if $w=1$. We now prove:

\begin{lemma}\label{formulas}
For any $w\neq 1$  in $P^\sigma$ and any reflection $\rho\in\WC$
such  that $\rho\preccurlyeq_\WA  w$ we  have $f(\rho  w)=f(w)-1$.
\end{lemma}
\begin{proof}
Recall that
$\sigma$  is induced  by the  bijection $a\mapsto  1-a$ of  $\bbZ$ so
that by Lemma \ref{reflexions de tilde C} a
reflection  in $\WC$ is either equal to  $(a,2kn+1-a)$ for some $a$ and $k$ in
$\bbZ$  or $(a,b)(1-a,1-b)$ for  some $a$ and  $b$ in $\bbZ$ with
$a+b\not\equiv 1 \text{ mod }2n$. Decompose
$w$ into a product $w=w_1w_2\dots w_r$ as in lemma \ref{formule pour l dans Ctilde}.
Let $\rho$ be a reflection
of $\WC$ which divides $w$ in $\WA$. By Proposition \ref{p_w isomorphism}, if
$\rho= (a,2kn+1-a)$  for some $a$, then
$a$ and $2kn+1-a$ have to be in the support of the same $w_i$ and if
$\rho=(a,b)(1-a,1-b)$, then $a$ and $b$ have to be in the support
of the same $w_i$ and $1-a$ and $1-b$ are then in the support of $\lexp\sigma w_i$.
Lemma \ref{formulas} is then a consequence of the following formulas, obtained by an immediate
computation. In each of these formulas the left hand side has
the form $\rho w_i$ with $w_i$ as above and
one checks easily that $f$ has value $f(w_i)-1$ on the
the right hand side.

\begin{multline*}
(a_1,2kn+1-a_1)(a_1,a_2,\ldots, a_h,2kn+1-a_1,2kn+1-a_2,\ldots, 2kn+1-a_h)=\\
\shoveright{(a_1,a_2,\ldots,a_h)(1-a_1,\ldots, 1-a_h)}\\
\shoveleft{(a_1,2kn+1-a_1)(a_1,\ldots,a_h)_{[1]}(1-a_1,\ldots,1-a_h)_{[-1]}=}\\
\shoveright{(a_1,\ldots,a_h,(2k+1)n+1-a_1,\ldots,(2k+1)n+1-a_h)}\\
\shoveleft{(a_1,b_1)(1-a_1,1-b_1)(a_1,\ldots,a_h,b_1,\ldots,b_k)
(1-a_1,\ldots,1-a_h,1-b_1,\ldots,1-b_k)=}\\
\shoveright{(a_1,\ldots,a_h)(b_1,\ldots,b_k)(1-a_1,\ldots,1-a_h)
(1-b_1,\ldots,1-b_k)}\\
\shoveleft{(a_1,b_1)(1-a_1,1-b_1)(a_1,\ldots,a_h,b_1,\ldots,b_k)_{[1]}
(1-a_1,\ldots,1-a_h,1-b_1,\ldots,1-b_k)_{[-1]}=}\\
\shoveright{(a_1,\ldots,a_h)(1-a_1,\ldots,1-a_h)(b_1,\ldots,b_k)_{[1]}
(1-b_1,\ldots,1-b_k)_{[-1]}}\\
\shoveleft{(a_1,b_1)(1-a_1,1-b_1)(a_1,\ldots,a_h,1-b_1,\ldots,1-b_k)_{[1]}
(1-a_1,\ldots,1-a_h,b_1,\ldots,b_k)_{[-1]}=}\\
{(a_1,\ldots,a_h,1-a_1,\ldots,1-a_h)
(b_1,\ldots,b_k,2n+1-b_1,\ldots,2n+1-b_k)}
\end{multline*}
\end{proof}

We can now finish the proof of Lemma \ref{formule pour l dans Ctilde} by induction on
$l_\WC(w)$.
If $w=\rho_1\rho_2\ldots\rho_l$ is a decomposition of $w$ as a product of reflections
in $\WC$, with $l=l_\WC(w)$, then $\rho_1\preccurlyeq_\WC w$, so that $\rho_1\preccurlyeq_\WA
w$ and by Lemma \ref{formulas} we have $f(\rho_1 w)=f(w)-1$.
Since $f(\rho_1 w)=l_\WC(\rho_1 w)$
by induction and $l_\WC(\rho_1 w)=l_\WC(w)-1$, we get the result.
\end{proof}

We now end the proof of Proposition \ref{divisibility restricts}.
First, since $l_\WC=f$ on $P^\sigma$,
Lemma \ref{formulas} shows that for any $w\in P^\sigma$ and any reflection $\rho\in\WC$
such that $\rho\preccurlyeq_\WA w$ we have $\rho\preccurlyeq_\WC w$.

Assume now that $w\preccurlyeq_\WA w'\preccurlyeq_\WA c$ with $w\neq 1$ and
$w,w'\in\WC$, so that $w,w'\in P^\sigma$;
choose a reflection  $\rho$ in $\WC$ such that $\rho\preccurlyeq_\WC w$, hence
$\rho\preccurlyeq_\WA w\preccurlyeq_\WA w'$ and $\rho w\preccurlyeq _\WA\rho w'$, so that in
particular $\rho\preccurlyeq_\WC w'$ by what we have just proved. 
We can assume by induction on $l_\WA(w)$ that $\rho w\preccurlyeq_\WC\rho w'$. Since
$l_\WC(\rho w)=l_\WC(w)-1$ and $l_\WC(\rho w')=l_\WC(w')-1$, we get $l_\WC(w)+l_\WC(w\inv
w')=1+l_\WC(\rho w)+l_\WC((\rho w\inv)( \rho w'))=1+l_\WC(\rho w')=l_\WC(w')$ so that
$w\preccurlyeq_\WC w'$.
\end{proof}

\begin{corollary}\label{presentation1}
\begin{enumerate}
\item
The germ $P^\sigma$ is in one-to-one correspondence with $\{w\in\WC\mid w\preccurlyeq_{\WC}c\}$
and this correspondence is compatible with divisibility.
\item\label{presentation duale}
The monoid $M(P^\sigma)$ has a presentation with set of generator $\uR$ in one-to-one correspondence
$\ur\mapsto r$ with the set $R$ of reflections $r\in\WC$ such that $r\preccurlyeq_\WC c$ and relations
$\ur_0\ur_1\ldots\ur_n=\ur'_0\ur'_1\ldots\ur'_n$ for any two $n$-tuples of reflections such that
$r_0\ldots r_n=r'_0\ldots r'_n=c$.
\end{enumerate}
\end{corollary}
Assertion (i) says that the germ $P^\sigma$ is obtained from $\WC$ and the balanced element
$c$ by the method of the ``generated group''.
\begin{proof}
Assertion (i) is precisely the content of Proposition \ref{divisibility
restricts}.

Assertion (ii) is always true for a monoid obtained by the method of the
generated group. We give a proof for the sake of completeness.
If $w\in\WC$ is such that $w\preccurlyeq_\WC c$, we denote by $\uw$ the corresponding element of the germ
$P^\sigma$. For such a $w$ there exist a sequence of reflections $r_0,\ldots,r_n$
and a non-negative integer $k\leq n$ such that $w=r_0r_1\ldots r_k$ and $c=r_0\ldots r_n$,
so that we have $\uw=\ur_0\ldots\ur_k$. A defining relation
of $M(P^\sigma)$ such as $\uw.\uw'=\underline{ww'}$ can be written
$\ur_0\ldots\ur_n=\uc$ if $w=r_0\ldots r_k$, $w'=r_{k+1}\ldots r_l$
and $c=ww'r_{l+1}\ldots r_n$, using the cancellability of the germ. This means that 
we can reduce the set of generators to $\uR$ and have the presentation given in (ii).
\end{proof}
\begin{definition}\label{dual monoid}
We call dual monoid of type $\tilde C_n$ the monoid $M(P^\sigma)$.
\end{definition}
\begin{remark}
Note that since all Coxeter elements are conjugate in $\WC$ changing the
Coxeter element in the presentation of the dual monoid
given in Corollary \ref{presentation duale} leads to isomorphic monoids.
\end{remark}
\section{Hurwitz action; presentations}
\label{hurwitz}
In this section we show that the Artin-Tits group of type $\tilde C_n$
is a Garside group: more precisely it is the group of fractions $G(P^\sigma)$
of the Garside monoid 
$M(P^\sigma)$. We also get simpler presentations for the monoid and the group
than in Corollary \ref{presentation duale}.
For this we study
the  Hurwitz action  on reduced  decomposition of  elements of $P^\sigma$.
Recall if $g$ is an element of some group $G$, the braid group with $n$ stands
acts on the set of decompositions of $g$ into a product of $n$ elements of $G$
in the following way: if $(g_1,\ldots,g_n)$ is such that $g=g_1g_2\ldots g_n$,
the  action  of  the  elementary  braid  $\bs_i$  maps  $(g_1,\ldots,g_n)$  onto
$(g_1,\ldots,g_{i-1},g_ig_{i+1}g_i\inv,g_i,g_{i+2},\ldots,g_n)$. This action is called
the Hurwitz action. Our first aim is to prove the following result where
we call ``reduced decomposition'' of an element $w\in\WC$ a decomposition of $w$
into a product of $l_\WC(w)$ reflections.
\begin{theorem}\label{hurwitz transitive} Let $c$ be a Coxeter element of $\WC$ and
$w$ be a divisor of $c$; then
the Hurwitz action is transitive on the set of reduced decompositions of
$w$.
\end{theorem}
\begin{proof}
Note  first that,  since a conjugate of  a reflection  is a reflection, the
Hurwitz action on reduced decompositions of an element $w$ is
well defined.

By  induction  on  $l_\WC(w)$  it  is  sufficient  to show that, starting with
some fixed reduced decomposition of $w$,  if $\rho$ is a
reflection  which divides $w$ we can get by the Hurwitz action a decomposition
of $w$ which begins with $\rho$.

Since by the Hurwitz action we can bring to the first place any reflection which
appears in a reduced decomposition of $w$,
it is enough to show
that  for any reflection $\rho$  which divides $w$ in  $\WC$ we can get by the
Hurwitz action a reduced decomposition of $w$ involving $\rho$.

Since all Coxeter elements are conjugate in $\WC$
(see \cite[\S 6, n$^\circ$1, Lemme 1]{bourbaki}) and since the Hurwitz action
is compatible with conjugation, we may assume that
$c=\sigma_0\sigma_2\sigma_4\ldots\sigma_1\sigma_3\sigma_5\ldots$ with the same
notation  as  in  Section  \ref{generated group}.  We view $\WC$ as $\WA^\sigma$ as
before.  We write $w=w_1\ldots w_r$ with $w_i$  as in \ref{formule pour l dans
Ctilde}.  We get a reduced decomposition of $w$ by concatenation of reduced
decompositions of the $w_i$.
A reflection  which divides  $w$ divides  some $w_i$, so
it is sufficient to show the result for $w=w_i$ which we now assume.

The element  $w$ is  either a  $\sigma$-stable finite  cycle, or  the product  of a
finite  cycle and its image by  $\sigma$ or a $\sigma$-stable pseudo-cycle. In
the first case $w$ can be written
$(a_1,\ldots,a_p,2kn+1-a_1,\ldots,2kn+1-a_p)$,   with $k\in\bbZ$,
all   $a_i$   odd  and
$a_1<a_2<\cdots<a_p<a_1+2n$,  so is a Coxeter element  of the Coxeter group of
type  $C_p$  generated  by  $(a_i,a_{i+1})(2kn+1-a_i,2kn+1-a_{i+1})$ for $i=1,
\ldots  p-1$ and  $(a_1,2kn+1-a_1)$. In  the second  case $w$ can be written
$(a_1,\ldots,a_r,\ldots,a_p)
(1-a_1,\ldots,1-a_r,\ldots,1-a_p)$ with $a_i$ odd for $i\leq r$ and $a_i$ even for $i>r$,
and the inequalities
$a_1<a_2<\cdots<a_r<a_1+2n$ and $a_{r+1}>\ldots>a_p>a_{r+1}-2n$,
so that $w$ is a Coxeter element of the Coxeter
group  of  type  $A_{p-1}$  generated  by $(a_i,a_{i+1})(1-a_i,1-a_{i+1})$ for
$i=1, \ldots p-1$. In the third case $w$ can be written
$(a_1,a_3\ldots,a_{2p-1})_{[1]}(a_2,a_4,\ldots,a_{2p})_{[-1]}$     with    the
$a_{2i-1}$  odd and $a_i+a_{2p+1-i}=2n+1$,  so that $w$
is a  Coxeter element of the
Coxeter group of type $\tilde C_p$ generated by $(a_1,a_{2p})$, the products
$(a_i,a_{i+1})(a_{2p-i},a_{2p+1-i})$ for $i=1,\ldots,p-1$ and $(a_p,a_{p+1})$.
Hence  $w$  is  a  Coxeter  element  of  a  Coxeter  subgroup $W$ of type $B$, $C$ or 
$\tilde C$; the reflections of $W$ are reflections of $\WC$  and a
reflection of $\WC$ divides $w$ in $\WC$ if and only if it divides $w$ in $W$.
A reduced decomposition of $w$ in $\WC$ is thus a reduced decomposition
of $w$ in that Coxeter subgroup.

Since  we know that the Hurwitz action is transitive on the reduced decomposition
of an element in groups of type $A$ or
$C$ (see \cite[Proposition 1.6.1]{bessis}), the only case which remains is the
case of a Coxeter element of a group of type $\tilde C$: we have only to study
the  case where $w$ is 
$c$ itself. We start with the reduced decomposition of $c$ given by
$(\sigma_0,\sigma_2,\sigma_4,\ldots,\sigma_1,\sigma_3,\sigma_5,\ldots)$.
We remark that if we delete $\sigma_0$ (resp.\ $\sigma_n$) from
this decomposition
we get a reduced decomposition of a Coxeter element of the Coxeter group
$W'$ of type $C_n$ generated by
$\sigma_1,\sigma_2,\ldots,\sigma_n$ (resp.\ of the Coxeter group
$W''$  of type $C_n$ generated by $\sigma_0$, $\sigma_1$,\dots,$\sigma_{n-1}$).
Using  the transitivity of the Hurwitz action for Coxeter groups of type $C_n$,
we know that by the Hurwitz action we can make appear any reflection of $W'$
or $W''$ in some reduced decomposition of $c$. 

We remark that if
$c=\rho_1\ldots\rho_{n+1}$  is a  reduced decomposition,  the Hurwitz orbit of
$(\rho_1,\ldots,\rho_{n+1})$  contains  $(\rho_2,\ldots,\rho_{n+1},c\inv\rho_1
c)$  so that by the  Hurwitz action for any $n\in\bbZ$ and any $i\in\{1,\ldots,n+1\}$
we can  get a reduced decomposition of $c$ involving  $c^n\rho_i c^{-n}$.
So we will get our result if we prove that any reflection appearing in a
reduced decomposition of $c$ is conjugate by some power of $c$ to a
reflection of $W'$ or of $W''$.

The action of $c$ by conjugation on permutations of
$\bbZ$  is induced by the translations
$2i\mapsto 2i-2$ and  $2i+1\mapsto 2i+3$ for all $i$.
The reflections of $W'$ are
$$\begin{cases}(a,b)(1-a,1-b)&
\text{ with }1\leq a<b\leq 2n\text{ and } a+b\neq 2n+1,\\
(a,2n+1-a)&\text{with }1\leq a\leq 2n.\end{cases}$$
The reflections of $W''$ are 
$$\begin{cases}(a,b)(1-a,1-b)&\text{ with }
1-n\leq a<b\leq n\text{ and }a+b\neq 1,
\\ (a,1-a)&\text{ with }1-n\leq a\leq n.\end{cases}$$

If $a$ and $b$ have same parity with $a<b\leq a+2n$,
we can conjugate $(a,b)(1-a)(1-b)\in P^\sigma$ by a power of $c$
to $(1,b+1-a)(0,a-b)$ which is a reflection of $W'$.
If $a$ and $b$ have different parities and $a+b\neq 1\text{ mod }2n$ we can 
conjugate $(a,b)(1-a,1-b)$ by a power of $c$ to 
$(\frac{a+b-1}2,\frac{a+b+1}2)(\frac{1-a-b}2,\frac{3-a-b}2)$ which
is a reflection of $W'$. Last of all we can conjugate 
$(a,2kn+1-a)$ where $k$ is arbitrary to $(kn,kn+1)$ which is equal either 
to $(n,n+1)$ which is in $W'$ or to $(0,1)$ which is in $W''$.
\end{proof}
\begin{remark}
\begin{enumerate}
\item
One can
conjecture that the Hurwitz action is transitive on the reduced decompositions
of a Coxeter element in any Coxeter group.
\item We have seen in the proof of Theorem \ref{hurwitz transitive} that a divisor of a Coxeter
element of $\WC$ is a Coxeter element of a Coxeter subgroup which is a direct product of groups
of type $A$, $C$ or $\tilde C$ (there can be at most one component of type $\tilde C$ and there
cannot be at the same time a component of type $C$ and a component of type $\tilde C$).
\end{enumerate}
\end{remark}

Using  Theorem  \ref{hurwitz  transitive}  we   can  simplify the relations
in the presentations  of
$M(P^\sigma)$  and $G(P^\sigma)$ given in Corollary \ref{presentation duale}.
As in Corollary \ref{presentation duale} and in its proof we denote by $\uw$ the element of $P^\sigma$
corresponding to an element $w\in\WC$ dividing $c$.
From this presentation we deduce one of our main result, which is that the Artin-Tits group
$G(\tilde C_n)$ is isomorphic to $G(P^\sigma)$. We use the notations $R$ and $\uR$ as
defined in Corollary \ref{presentation1}.

\begin{theorem}\label{presentation}
\begin{enumerate}
\item
The monoid $M(P^\sigma)$ has the following monoid presentation by generators and relations:
$$M(P^\sigma)=<\uR\mid \ur.\ut=\underline{rtr}.\ur \text{ if }r,t\in R \text{ and }
rt\preccurlyeq_\WC c>^+$$
\item
The Artin-Tits group of type $\tilde C_n$ is isomorphic to $G(P^\sigma)$; in 
particular it is the Garside group of the Garside monoid $M(P^\sigma)$ and has
the group presentation
$G(\tilde C_n)=<\uR\mid \ur.\ut=\underline{rtr}.\ur \text{ if }r,t\in R \text{ and }
rt\preccurlyeq_\WC c>$.
\end{enumerate}
\end{theorem}
\begin{proof}
By Corollary 
\ref{presentation duale} the monoid $M(P^\sigma)$ has a presentation with set of generators 
$\uR$ and relations given by
the equalities between the lifts of any two reduced decompositions of
$c$. We have to prove
that one can pass from a reduced decomposition of $c$
to another one by applying only relations  of the form
$ r.t=(rtr).r$. But this is precisely the transitivity of the Hurwitz action, whence (i).

Since all the simple reflections of $\WC$ are in $R$, we have elements
$\usigma_i\in\uR$ for $i=0,\ldots,n$,
where as in Section \ref{generated group} the simple reflections of $\WC=\WA^\sigma$
are $\sigma_i=s_is_{2n-i}$ 
for $i=2,\ldots,n-1$, $\sigma_n=s_n$ and $\sigma_0=s_{2n}$.
The natural morphism $M(P^\sigma)\to M(P)^\sigma$ extends to a group morphism $f:G(P^\sigma)\to
G(P)^\sigma$. We know that $G(P)$ is the group $\Atilde$ so that $G(P)^\sigma$ is the group 
$G(\tilde C_n)$ by Theorem \ref{G(tilde C)=C_Atilde(sigma)}. The morphism $f$ maps $\usigma_i$
to $\bsigma_i$ for $i=2,\ldots,n-1$ and maps $\usigma_n$ to $\bsigma_n$ and $\usigma_0$
to $\bsigma_0$.

We claim that $\usigma_0$,\dots, $\usigma_n$ satisfy the braid relations of type $\tilde C_n$. 
Indeed, for any pair of distinct simple reflections $\sigma_i$ and $\sigma_j$,
we have $\sigma_i\sigma_j\preccurlyeq_\WC c$, up to swapping $i$ and $j$; 
if $\sigma_i$ and $\sigma_j$, commute, then
$\sigma_i\sigma_j\preccurlyeq_\WC c$ so that $\usigma_i\usigma_j=\usigma_j\usigma_i$
in $P^\sigma$; if $\sigma_i$ and $\sigma_j$ satisfy a braid relation of length 3 we have
$\usigma_j.\usigma_i.\usigma_j=\usigma_j.\underline{\sigma_i\sigma_j\sigma_i}.\usigma_i=
\underline{\sigma_j\sigma_i\sigma_j\sigma_i\sigma_j}.\usigma_j.\usigma_i=
\usigma_i.\usigma_j.\usigma_i$; if $\sigma_i$ and $\sigma_j$ satisfy a braid relation of length
4 we have
$\usigma_i.\usigma_j.\usigma_i.\usigma_j=\usigma_i.\usigma_j.\underline{\sigma_i\sigma_j\sigma_i}.
\usigma_i=\usigma_i.\underline{\sigma_j\sigma_i\sigma_j\sigma_i\sigma_j}.\usigma_j.\usigma_i=
\underline{\sigma_i\sigma_j\sigma_i\sigma_j\sigma_i\sigma_j\sigma_i}.\usigma_i.\usigma_j.\usigma_i=
\usigma_j.\usigma_i.\usigma_j.\usigma_i$. Hence there exists a morphism $g:G(\tilde C_n)\to
G(P^\sigma)$ which maps $\bsigma_i$ to $\usigma_i$ for $i=0,\ldots,n$. We have then $f\circ g=\Id$
since $G(\tilde C_n)$ is generated by $\{\bsigma_i,i=0\ldots,n\}$.

On the other hand the transitivity of the Hurwitz action and the presentation of $G(P^\sigma)$
imply that one can express any $\ur$ where $r\preccurlyeq_\WC c$ as a conjugate of some
$\usigma_i$ by a product of elements of the form $\usigma_j^{\pm1}$, so that $G(P^\sigma)$ is
generated by $\{\usigma_i, i=0,\ldots,n\}$. This implies that $g$ is surjective, so that $f$ and
$g$ are isomorphisms.
\end{proof}

\section{Some consequences}
Among consequences of our results we can recover the known fact that the center of $B(\tilde C_n)$ is
trivial. The Garside structure gives also a solution to the word problem in this group.
As an illustration of the use of a Garside structure we compute the centralizer of powers of a Coxeter
element in an Artin-Tits group of type $\tilde C_n$.
\begin{proposition}
With the notation of section \ref{Artin-Tits de type tilde C}, let
$\bc=\bsigma_0\bsigma_1\ldots\bsigma_n$; then for any $h\in\bbZ-\{0\}$, the centralizer of $\bc^h$ in
the Artin-Tits group $G(\tilde C_n)$ is isomorphic to the Artin-Tits group of type $C_{\gcd(h,n)}$.
\end{proposition}
\begin{proof}
It is a general result (see \eg, \cite[Proposition 2.26]{BDM})
that the monoid of fixed points $M^\varphi$ under an automorphism
$\varphi$  of  a  Garside  monoid  $M$ fixing the Garside
element $\Delta$ 
is a Garside  monoid with same Garside element (this can also be seen as a particular
case of Corollary \ref{points fixes}): if $P$ is the set of
divisors of $\Delta$ in $M$ then the set of divisors of $\Delta$ in $M^\varphi$ is
$P^\varphi$. In this situation
the group of fixed points of $\varphi$ in the group $G$ of fractions of $M$
is the group of fractions of $M^\varphi$.
We apply this to $M=M(P^\sigma)$ and $G=G(P^\sigma)$, taking for $\varphi$ the
conjugation  by  $\bc^h$.  As  in  section  \ref{fixed points} we identify the
elements  of  $P^\sigma$  with  the  elements  of  $\WC$  corresponding to the
$2n$-periodic,  non-crossing $\sigma$-stable  partitions of  $X\cup\Xi$, where
$X=1+2\bbZ$  and $\Xi=2\bbZ$. Recall that $\bc$ is the permutation
given by $x\mapsto x+2$ for $x\in X$ and $\xi\mapsto \xi-2$ for $\xi\in\Xi$.

Let  $p$  be  a
$\sigma$-stable non-crossing partition stable by the action of $\bc^h$. If $A$
is  a part of $p$ contained in $X$ or  in $\Xi$, then $A+2kh$ is also a
part of $p$ for all $k\in\bbZ$.  If $A$ is a part of $p$
which intersects both $X$ and $\Xi$ and is not $\bc^h$-stable, then $A$
crosses its image under the action of $\bc^h$, which is impossible since two
distinct  parts of $p$ are  non-crossing. Hence any part $A$ of $p$ has to  be
invariant  by the  action of
$\bc^h$  which is equivalent  to $A=A+2kh$ for  all $k\in\bbZ$. Hence a part of a
partition $p_w$
associated  to  an  element  $w\in (P^\sigma)^{\bc^h}$ (see Proposition \ref{p_w
isomorphism}) is either an
infinite part both $2n$ and $2h$ periodic (and there is at most one such part in $p_w$)
or a finite part $A$ contained in $X$ or in $\Xi$ and
such that $A+2kh$ and $A+2kn$ are parts of $p_w$ for any $k$. 
Conversely if all the parts of a $\sigma$-stable non-crossing partition $p_w$
are of one of these form then $w$ is centralized by $\bc^h$.

From the above study we see that the centralizer of $\bc^h$ is equal to the centralizer
of $\bc^{\gcd(n,h)}$. Hence to prove the proposition we are reduced to the case where $h$
divides $n$.

We first study the case $h=n$.
\begin{lemma}
The centralizer of $\bc^n$ is an Artin-Tits group of type $C_n$.
\end{lemma}
\begin{proof}
From the above discussion we see that
an element $w\in (P^\sigma)^{\bc^n}$ is a product of pairwise commuting
and non-crossing elements of the form either
$(2a_1+1,2a_2+1,\ldots,2a_r+1)(-2a_1,-2a_2,\ldots,-2a_s,\ldots,-2a_r)$
with $a_1<a_2<\ldots<a_r<a_1+n$, or
$(2a_1+1,2a_2+1,\ldots,2a_r+1)_{[1]}(-2a_1,-2a_2,\ldots,-2a_r)_{[-1]}$
with $0< a_1<a_2<\ldots<a_r\leq n$.
Hence all elements of $(P^\sigma)^{\bc^n}$ can be written $x.\lexp\sigma x$ where $x$ is a
permutation of the set $X$ of 
odd integers and $\lexp\sigma x$ is a permutation of the set $\Xi$ of even
integers. Since a permutation of $X$ and a permutation of $\Xi$ are always
non-crossing,
two such elements $x.\lexp\sigma x$ and $y.\lexp\sigma y$ in $(P^\sigma)^{\bc^n}$
are non-crossing if and only if $x$ and $y$ are non-crossing.
In \cite[Proof of 5.10]{Atilde} we have described the germ for the dual presentation
of an Artin-Tits group of type $C_n$ as the set of $n$-periodic, non-crossing
permutations of $\bbZ$. From this description we see that
under the bijection $2k+1\mapsto k$  from $X$ to $\bbZ$, the map $x.\lexp\sigma x\mapsto x$
defines an isomorphism of germs from $(P^\sigma)^{\bc^n}$ to the dual germ
for the Artin-Tits group of type $C_n$.
\end{proof}
We deduce the proposition using the fact that the centralizer of $\bc^h$ for a divisor $h$ of $n$ is the
centralizer of $\bc^h$ in the centralizer of $\bc^n$ so that it is the centralizer of $\bc^h$ in an Artin
group of type $C_n$, which is known to be an Artin-Tits group of type $C_h$ (see \cite{BDM}).
\end{proof}
\bibliographystyle{amsplain}
\bibliography{Cntilde}
\end{document}